\documentclass[graybox]{svmult}

\usepackage{type1cm}        % activate if the above 3 fonts are
\usepackage{makeidx}         % allows index generation
\usepackage{graphicx}        % standard LaTeX graphics tool
\usepackage{multicol}        % used for the two-column index
\usepackage[bottom]{footmisc}% places footnotes at page bottom

\usepackage{newtxtext}       %
\usepackage[varvw]{newtxmath}       % selects Times Roman as basic font

\usepackage{xcolor}
\usepackage{comment}
\usepackage{tikz}

\newcommand{\fonction}[5]{
    \begin{array}{l@{\,}ccc}
      #1\colon
      & #2
      & \longrightarrow
      & #3
      \\
      & \displaystyle{#4}
      & \longmapsto
      & \displaystyle{#5}
    \end{array}}

\newcommand{\limit}[2]{{\ \underset{#1 \to #2}{\longrightarrow} \ }}

\newcommand{\R}{{\mathbb{R}}}

\newcommand{\Ck}[1]{\mathbf{C}^{#1}}
\newcommand{\Cc}[1]{\mathbf{C}_c^{#1}}
\renewcommand{\L}[1]{{\mathbf{L}^#1}}
\newcommand{\Lloc}[1]{{\mathbf{L}_{\mathbf{loc}}^{#1}}}

\newcommand{\p}{\partial}
\newcommand{\sgn}{\mathop{\rm sgn}}
\renewcommand{\d}[1]{\mathinner{\mathrm{d}{#1}}}
\newcommand{\flow}{\mathcal{F}}

\newcommand{\red}[1]{\textcolor{red}{#1}}
\newcommand{\ds}[1]{\displaystyle{#1}}

\makeindex             % used for the subject index
\begin{document}

\title*{Peculiarities of Space Dependent Conservation Laws:
    Inverse Design and Asymptotics} \titlerunning{Inverse Design
    and Asymptotics in Space Dependent Conservation Laws}
% Use \titlerunning{Short Title} for an abbreviated version of your
% contribution title if the original one is too long
\author{Rinaldo M. Colombo, Vincent Perrollaz and Abraham Sylla}
\authorrunning{R. M. Colombo, V. Perrollaz and A. Sylla}
% Use \authorrunning{Short Title} for an abbreviated version of your
% contribution title if the original one is too long
\institute{Rinaldo M.~Colombo \at Unit\`a INdAM \& Dipartimento di Ingegneria
dell'Informazione, Universit\`a di Brescia\\
  \email{rinaldo.colombo@unibs.it} \and Vincent Perrollaz \at Institut
  Denis Poisson, Universit\'e de Tours, CNRS UMR 7013, Universit\'e
  d'Orl\'eans\\ \email{Vincent.Perrollaz@univ-tours.fr} \and
  Abraham Sylla \at  Dipartimento di Matematica e Applicazioni,
  Università di Milano-Bicocca \\
  \email{abraham.sylla@unimib.it}}
%
% Use the package "url.sty" to avoid problems with special characters

\maketitle

\begin{abstract}
  {Recently, results regarding the Inverse Design problem for
    Conservation Laws and Hamilton-Jacobi equations with $x$-dependent
    convex fluxes were obtained in~\cite{CPS2022, CPS2023}. More
    precisely, characterizations of attainable sets and the set of
    initial data evolving at a prescribed time into a prescribed
    profile were obtained. Here, we present an explicit example that
    underlines deep differences between the $x$-dependent and
    $x$-independent cases. Moreover, we add a detailed analysis of the
    time asymptotic solution of this example, again underlining
    differences with the $x$-independent case.}
\end{abstract}

\section{Introduction}
\label{sec:Intro}

Consider the scalar one dimensional conservation law
\begin{equation}
  \label{eq:cl}\tag{\textbf{CL}}
  \left\{
    \begin{array}{l@{\qquad}r@{\,}c@{\,}l}
      \p_t u + \p_{x} (H(x,  u)) = 0 & (t,x)
      & \in & \left]0,+\infty\right[ \times \R \\
      u(0,x) = u_o(x) & x & \in & \R\,.
    \end{array}
  \right.
\end{equation}
Denote by
\begin{equation}
  \label{eq:semigroup}
  S^{CL} \colon \R^+ \times \L{\infty}(\R; \R) \to \L{\infty}(\R; \R)
\end{equation}
the semigroup whose orbits are the entropy solutions
to~\eqref{eq:cl}. For any positive $T$ and for any assigned profile
$w \in \L{\infty}(\R; \R)$ the inverse design is the set of initial
data evolving into this profile at time $T$, i.e.,
\begin{equation}
  \label{eq:InverseDesign}
  \begin{array}{rcl}
    I_T^{CL}(w)
    & \coloneqq
    & \left\{u_o \in \L{\infty}(\R; \R) \colon S_T^{CL}u_o = w \right\}.
  \end{array}
\end{equation}
In the homogeneous case, a general characterization of $I^{CL}_T (w)$
is given in~\cite{CP2020}. Other more specific results in this setting
are~\cite{LZ2021}, devoted to Burgers' equation; \cite{MR1616586},
specific to the attainable set for boundary value problems arising in
the modeling of vehicular traffic.  However, the homogeneous case
significantly differs from the present non-homogeneous one and less
results in the literature are available. A first step in this
direction, limited to the study of the attainable set,
is~\cite{MR4188826}, see also the related
preprint~\cite{AdimurthiGhoshal2020}, where $H$ in~\eqref{eq:cl}
consists of an expression for $x>0$ and another expression for $x<0$.
A specific inverse design problem related to conservation laws, e.g.,
the minimization of a sonic boom, is considered in~\cite{MR3643881},
while~\cite{MR4503821} is devoted also to Hamilton-Jacobi equation.

In~\cite{CPS2023}, the characterizations obtained in~\cite{CP2020} is
extended to the non-homogeneous case. The analytic techniques
developed in \cite{CPS2023} take advantage of the connection
$u = \partial_x U$ between \eqref{eq:cl} and the Hamilton-Jacobi
equation
\begin{equation}
  \label{eq:hj}\tag{\textbf{HJ}}
  \left\{
    \begin{array}{l@{\qquad}r@{\,}c@{\,}l}
      \p_t U + H(x,  \p_{x}U) = 0 & (t,x) & \in & \left]0,+\infty\right[ \times \R \\
      U(0,x) = U_o(x) & x & \in & \R\,.
    \end{array}
  \right.
\end{equation}
We know, on the basis of~\cite{CPS2022}, that both Cauchy problems for
\eqref{eq:cl} and \eqref{eq:hj} are globally well posed under the same
set of assumptions, namely:
\begin{align}
  \label{eq:smoothness}\tag{\textbf{C3}}
  \mathbf{Smoothness}:
  & \qquad H = H(x, u) \in \Ck{3}(\R^2, \R).
  \\[6pt]
  \label{eq:heterogeneity}\tag{\textbf{CNH}}
  \mathbf{Compact~NonHomogeneity:}
  & \qquad
    \begin{array}{@{}l}
      \exists X > 0, \; \forall (x,u) \in \R^2
      \\
      |x| > X \implies \p_x H (x,u) = 0.
    \end{array}
  \\[6pt]
  \label{eq:convexity}\tag{\textbf{CVX}}
  \mathbf{Strong~Convexity:}
  & \qquad
    \begin{array}{@{}l}
      \forall x \in \R, \; u \mapsto \p_u H(x, u) \mbox{ is an}
      \\
      \mbox{increasing } \Ck{1} \mbox{-diffeomorphism}
      \\
      \mbox{of } \R \mbox{ onto itself.}
    \end{array}
\end{align}
Assumption~\eqref{eq:convexity} implies that $H$ is strictly convex
with respect to the second variable. Clearly, the strongly concave
case is entirely analogous. Let us also mention that for
well-posedness, \eqref{eq:convexity} is relaxed to a uniform
coercivity assumption coupled to a genuine nonlinearity assumption
in~\cite[\textbf{(UC)-(WGNL)}]{CPS2022}. Rather than tackling directly
the characterization of the inverse design for \eqref{eq:cl},
\cite{CPS2023} deals first with~\eqref{eq:hj} and then uses the
correspondence to get back to~\eqref{eq:cl}.

Below, we construct an explicit example showing that when $H$ depends
on $x$, the inverse design $I^{CL}_T (w)$ may have properties, in a
sense, opposite to those that hold in the homogeneous case, according
to~\cite{CP2020}. Indeed, in the $x$-independent case, the presence of
a shock in $w$ implies that $I^{CL}_T (w)$ is a cone with infinite
dimensional extremal faces. In our $x$-dependent example, in spite of
$w$ displaying a shock, $I^{CL}_T (w)$ is a singleton. On the
contrary, in the homogeneous case, $I^{CL}_T (w)$ is a singleton if
and only if $w$ is continuous.

Furthermore, in this example, the time asymptotic solution is neither
constant nor a rescaling of a solution to a Riemann Problem. On the
contrary, as $t \to +\infty$, the solution converges locally uniformly
to a compactly supported profile displaying a single singularity,
which is a stationary entropic shock wave. This suggests various questions
about the characterization of time asymptotic profiles to non homogeneous
conservation laws and, in particular, about their stability properties.

\section{Notations and Results}

Recall the classical definition of entropy
solution~\cite[Definition~1]{Kruzhkov1970}, as tweaked
in~\cite{CPS2022}.

\begin{definition}
  \label{def:EntropySolution}
  Fix $u_o \in \L{\infty}(\R;\R)$. A bounded function
  $u \in \L{\infty}(\R^+ \times \R; \R)$ is a \emph{solution}
  to~\eqref{eq:cl} if for all test functions
  $\varphi \in \Cc{\infty}(\R^+ \times \R; \R^+)$ and for all scalar
  $k \in \R$:
  \begin{equation}
    \label{eq:EntropyIneq}
    \begin{array}{rcc}
      \ds{\int_0^{+\infty} \int_\R |u(t, x) - k| \, \p_t \varphi(t,x) \d{x}\d{t}}
      & & \\[5pt]
      \ds{+ \int_0^{+\infty} \int_\R \sgn \left(u(t, x) - k\right) \;
      \left(H \left(x,u(t, x)\right) - H (x,k)\right)
      \, \p_x \varphi (t,x) \d{x} \d{t}} & & \\[5pt]
      \ds{- \int_0^{+\infty}
      \int_\R \sgn\left(u(t, x) - k\right) \, \p_x H (x,k) \, \varphi (t,x)
      \d{x} \d{t}} & & \\[5pt]
      \ds{+ \int_\R |u_o(x) - k| \, \varphi(0,x) \d{x}} & \geq & 0.
    \end{array}
  \end{equation}
\end{definition}
Definition~\ref{def:EntropySolution}, taken from
by~\cite[Definition~2.1]{CPS2022} is, apparently, weaker than the
classical {Kru\v zkov} definition since it does not require the
\emph{``trace at $0$ condition''}~\cite[Formula~(2.2)]
{Kruzhkov1970}. Nevertheless, under Assumption~\eqref{eq:smoothness},
Definition~\ref{def:EntropySolution} ensures uniqueness and uniform
$\Lloc1$-continuity in time of the solution, as proved
in~\cite[Theorem~2.6]{CPS2022}.

As usual, in connection with~\eqref{eq:cl} we use the system of
ordinary differential equations
\begin{equation}
  \label{eq:HS}\tag{\textbf{HS}}
  \left\{
    \begin{aligned}
      \dot q & = \partial_p H(q, p) \\
      \dot p & = - \partial_q H(q, p),
    \end{aligned}
  \right.
\end{equation}
which we consider equipped either with initial or with final
conditions. For all $w \in \L{\infty}(\R; \R)$ such that
$I_T^{CL}(w) \neq \emptyset$, define
\begin{equation}
  \label{eq:foot}
  \fonction{\pi_w}{\R}{\R}{x}{q(0)}
  \mbox{ where }
  (q,p) \mbox{ solves~\eqref{eq:HS}}
  % \mbox{ with datum }
  \mbox{ with }
  \left\{
    \begin{array}{r@{\,}c@{\,}l}
      q (T) & = & x \\
      p (T) & = & w (x) \,.
    \end{array}
  \right.
\end{equation}
We also introduce the set
\begin{displaymath}
  %\label{eq:hamiltonian_rays}
  \mathcal{R}_T \coloneqq
  \left\{
    q \in \mathbf{C}^{1}([0,T];\R)  \colon
    \exists \, p \in \mathbf{C}^1([0,T]; \R) \text{ such that } (q, p)
    \mbox{ solves~{\eqref{eq:HS}}}
  \right\} .
\end{displaymath}
whose elements we call Hamiltonian rays. The map $\pi_w$ assigns to $x \in \R$ the intersection of the minimal
backward characteristics emanating from $(T, x)$,
see~\cite[Definition~3.1, Theorems~3.2-3.3]{Dafermos1977}, with the
axis $t=0$. Remark that in the $x$-independent case, all Hamiltonian
rays are straight lines, as also any extremal characteristics, a key
simplification exploited in~\cite[Formula~(2.3)]{CP2020}.

\begin{theorem}
  \label{th:InverseDesign}
  {\rm{\cite[Corollary 3.5]{CPS2023}}} Let $H$
  satisfy~\eqref{eq:smoothness}, \eqref{eq:heterogeneity}
  and~\eqref{eq:convexity}. Fix $T>0$ and $w \in \L\infty(\R; \R)$.
  Then, $I_T^{CL}(w) \neq \emptyset$ if and only if $\pi_w$ admits a
  non decreasing representative. In this case, $I_T^{CL}(w)$ is a
  closed convex cone with a unique vertex $u_o^*$.
\end{theorem}

A characterization of the attainable set for~\eqref{eq:cl} can be
obtained by that for~\eqref{eq:hj} provided
in~\cite[Theorem~3.2]{CPS2023}.  The latter theorem extends to the
$x$-dependent case some of the properties known to hold in the
$x$-independent case, see~\cite{CP2020}. However, the extension to the
$x$-dependent case can not be merely reduced to the rise of technical
difficulties. Indeed, some properties are irremediably lost and new
phenomena arise.

\begin{figure}[!h]
  \begin{center}
    \includegraphics[scale=0.33]{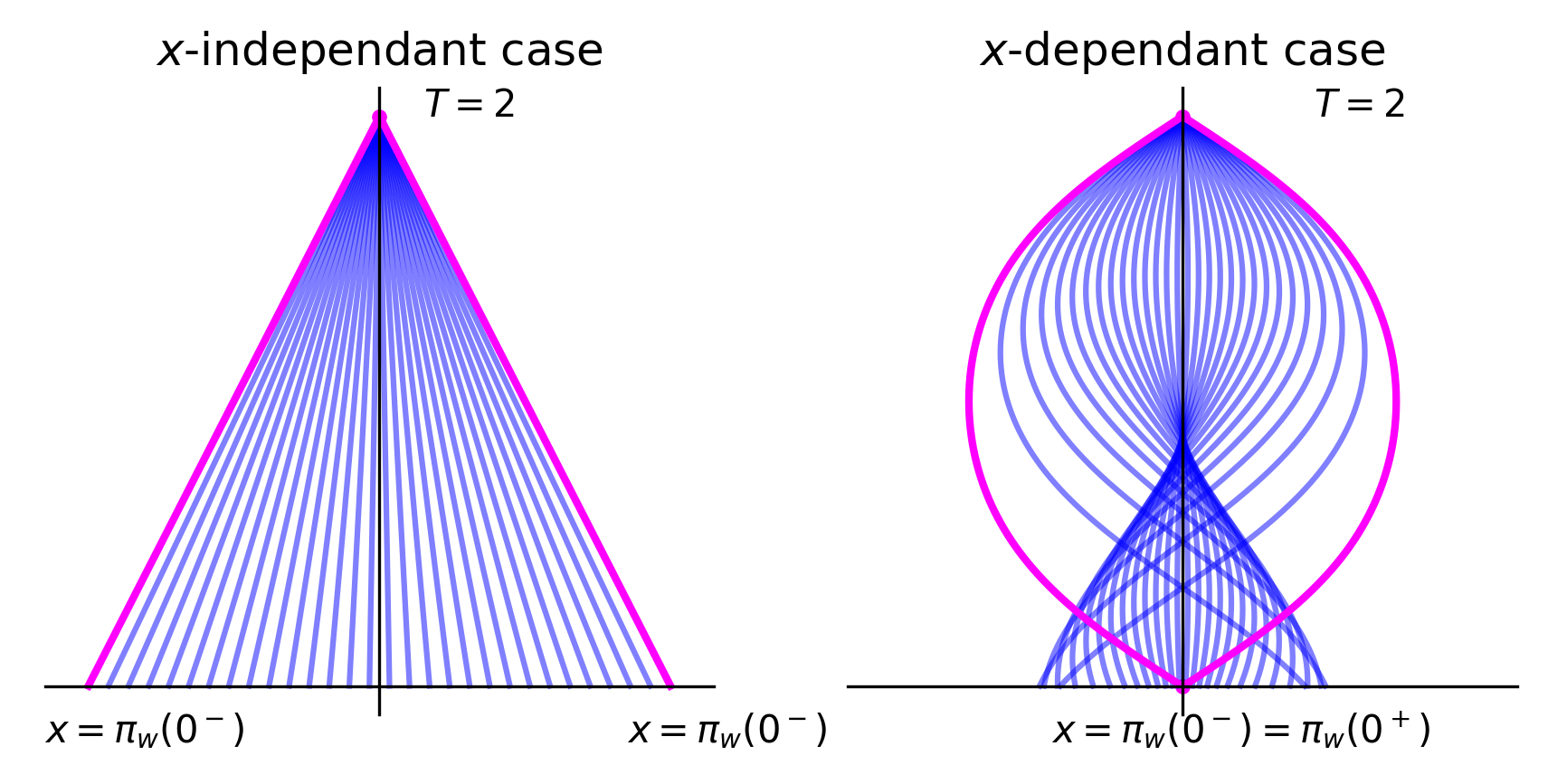}
  \end{center}
  \vspace*{-5mm}
  \caption{Left, in the $x$-independent case, the Hamiltonian rays
    fill the non uniqueness gap. Right, in the $x$-dependent case
    defined by the Hamiltonian~\eqref{eq:hamiltonian}, extremal
    characteristics still do not intersect, but Hamiltonian rays do
    and may well exit the non uniqueness gap or also intersect.}
  \label{fig:Extremals}
\end{figure}

The most apparent difference between the two situations is described
in Figure~\ref{fig:Extremals}, with reference to extremal backward
generalized characteristics (thicker curves), whose behaviors in the
two cases are quite different. In the $x$-independent case, extremal
backward characteristics define a \emph{non uniqueness gap}. On the
contrary, in the $x$-dependent case, extremal backward characteristics
may well intersect at the initial time, so that the non uniqueness gap
disappears. Furthermore, in the $x$-independent case, an isentropic
solution, see~\cite[Theorem~3.1]{CPS2022}, can be constructed filling
the non uniqueness gap with Hamiltonian rays, that is solutions to
\eqref{eq:HS}, emanating from $q(T) = 0$,
$p (T) = (1-\lambda) w(0-) + \lambda w(0+)$, for $\lambda \in
[0,1]$. On the contrary, the same procedure might not work in the
$x$-dependent case, see also~\cite{MR1722801} for a multi-d study. The
numerical integrations in Figure~\ref{fig:Extremals}, right, referred
to~\eqref{eq:HS} with Hamiltonian~\eqref{eq:hamiltonian}, show that
extremal backward characteristics still do not intersect in
$\mathopen]0,T\mathclose[ \times \R$, but the intermediate Hamiltonian
rays may well cross each other and even exit the region bounded by the
extremal characteristics.

When $H$ does not depend on $x$, $u_o^*$ defined in
Theorem~\ref{th:InverseDesign} is characterized by~\cite[(G2) in
Proposition~5.2]{CP2020}. Then, \cite[(R1) in Lemma~7.2]{CP2020}
ensures not only that $u_o^*$ is one sided Lipschitz continuous,
but also that the solution $\tilde u$ to~\eqref{eq:cl} with initial
datum $u_o^*$ evolving into $w$ is Lipschitz continuous on any
compact subset of $\mathopen]0, T \mathclose[ \times \R$. Thus,
$\tilde u$ satisfies the inequality in
Definition~\ref{def:EntropySolution} with an equality, \textit{i.e.},
it is an \emph{isentropic} solution and also reversible in time. All
this is no longer true in the heterogeneous case, as highlighted by
the following theorem.

\begin{figure}[!htp]
  \begin{center}
    \includegraphics[scale=0.42]{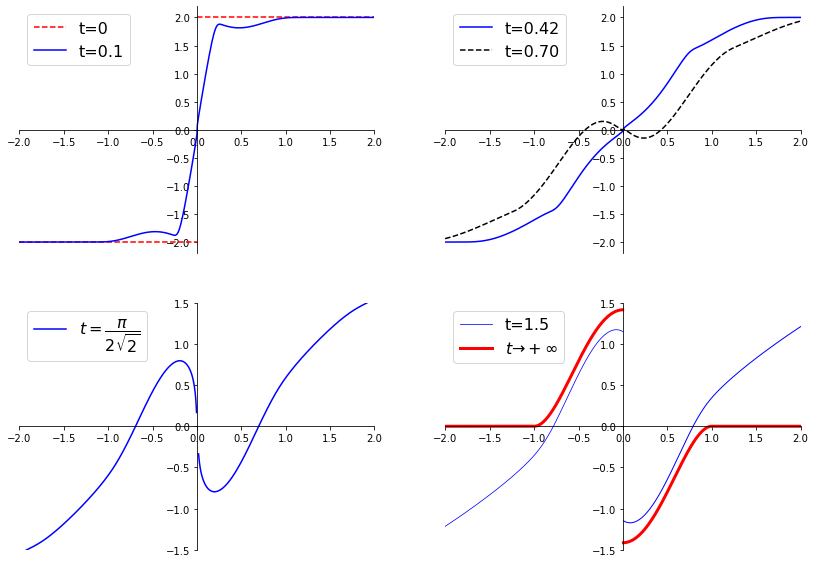}
  \end{center}
  % \vspace*{-5mm}
  \caption{Numerically computed solution to~\eqref{eq:cl} with $H$ as
    in~\eqref{eq:hamiltonian} and initial datum~\eqref{eq:initial}, as
    described in Theorem~\ref{th:example}. The Riemann initial
    datum~\eqref{eq:initial} first evolves similarly to a rarefaction
    (top left); then gets non monotone (top right); at time
    $t = \pi / (2\sqrt2)$ a shock is formed (bottom left) and in
    subsequent times the size of the shock grows (bottom right). The
    latter graph also portrays the asymptotic profile~\eqref{eq:2}.}
  \label{fig:NumSim}
\end{figure}

\begin{figure}[!h]
  \begin{center}
    \includegraphics[scale=0.45]{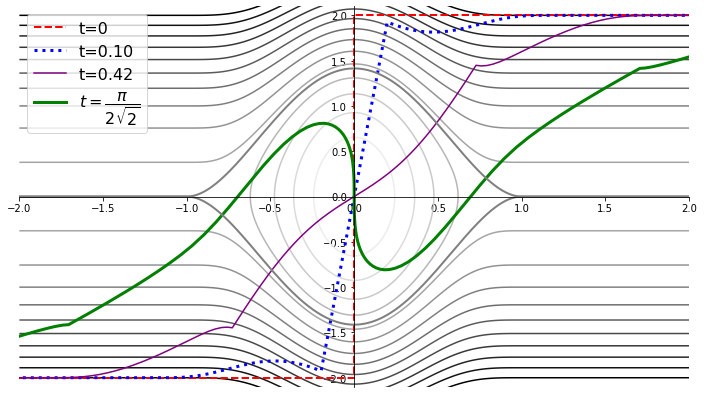}
  \end{center}
  \caption{Illustration of Proposition~\ref{prop:ce6}: evolution in
    time of the solution to~\eqref{eq:cl} in connection with the
    orbits of~\eqref{eq:HS}--\eqref{eq:hamiltonian}.}
  \label{fig:EvolFlow}
\end{figure}

\begin{theorem}
  \label{th:example}
  {\rm{\cite[Theorem 4.1]{CPS2023}}} Define
  \begin{equation}
    \label{eq:hamiltonian}
    H(x, p) \coloneqq  \frac{p^2}{2} + g(x)
    \quad \mbox{ where } \quad
    g(x) \coloneqq
    \left\{
      \begin{array}{cl}
        1 - (1 - x^2)^4 & \text{if} \;\; |x| \leq 1 \\
        1 & \text{otherwise.}
      \end{array}
    \right.
  \end{equation}
  Then, \eqref{eq:smoothness}, \eqref{eq:heterogeneity} and~\eqref{eq:convexity} hold. Moreover, there exists $w \in \L\infty (\R; \R)$ such that \\[5pt]
  $\mathbf{(i)}$ For all $T > 0$, $I_T^{CL}(w)$ is a singleton. \\[5pt]
  $\mathbf{(ii)}$ For all $T > \pi / (2\sqrt2)$ the solution $u$
  to~\eqref{eq:cl} with $u(0)$ being the vertex of $I^{CL}_T (w)$
  displays an entropic shock sited at $x = 0$, whose size grows in
  time, that arises at time $t = \pi / (2\sqrt2)$.
\end{theorem}
Elements of proof are presented in Section~\ref{sec:proof}. As we
mentioned before, the contrary to \textbf{(ii)} holds when $H$ does
not depend on $x$. In the homogeneous case, the contrary to
\textbf{(i)} also holds true, see~\cite[(G1) in
Proposition~5.2]{CP2020}.  The evolution of the numerical solution to
\eqref{eq:cl}-\eqref{eq:hamiltonian} with initial data in
$I_T^{CL}(w)$ (which is a singleton) computed with a standard finite
volume scheme, is represented in Figure~\ref{fig:NumSim}, see also
Figure \ref{fig:EvolFlow}.
% Remark, and this is intrinsic to the heterogeneous case, that the
% initial rarefaction profile evolves into a shock wave.

Therefore, there exist an $x$-dependent Hamiltonian $H$, a profile $w$
and a time $T > 0$ such that $I^{CL}_T(w) \neq \emptyset$ but in any
solution evolving from an initial datum in $I^{CL}_T(w)$, shocks arise
before time $T$, so that no reversible solution is possible, see
Figure~\ref{fig:EvolFlow}. In other words, the profile $w$ can be
reached exclusively producing a sufficient amount of entropy and no
isentropic solution evolves into $w$. Each of these facts requires $H$
to depend on $x$ and is false in an $x$-independent setting. As a
consequence, no direct definition of $u_o^*$ is available, as it
was in the $x$-independent case, explaining why in~\cite{CPS2023} the
authors had to resort to~\eqref{eq:hj} for its construction.

A further difference between the homogeneous and the non homogeneous
cases is in the long time behavior.

\begin{theorem}
  \label{th:asymptotics}
  Let $H$ be as in~\eqref{eq:hamiltonian}, choose
  \begin{equation}
    \label{eq:initial}
    u_o(x) \coloneqq
    \left\{
      \begin{array}{cl}
        -2 & \text{if } x < 0 \\
        2 & \text{if } x > 0\,,
      \end{array}
    \right.
  \end{equation}
  and consider the corresponding entropy solution $u$
  to~\eqref{eq:cl}. Then, for all $x \in \R \backslash \{0\}$,
  \begin{equation}
    \label{eq:2}
    \lim_{t\to+\infty} u (t,x)
    =
    -{\mathop {\rm sgn}}\, x \; \sqrt{2 (1-g (x))}
    =
    \left\{
      \begin{array}{l@{\qquad}r@{\,}c@{\,}l}
        - \sqrt2 \, {\mathop {\rm sgn}}\, x \; (1-x^2)^2
        & |x|
        & \leq
        & 1 \,,
        \\
        0
        & |x|
        & >
        & 1 \,.
      \end{array}
    \right.
  \end{equation}

  % ****

  % $t \mapsto u(t, x)$ admits a finite limit $u_{\infty}(x)$ as
  % $t \to \infty$. Moreover, \\[5pt]
  % (i) $u_\infty$ is an odd function.\\[5pt]
  % (ii) $u_\infty \in \L{\infty}(\R; [-\sqrt{2},
  % \sqrt{2}])$. \\[5pt]
  % (iii) $u_\infty$ is nondecreasing on $]0,+\infty[$.
\end{theorem}

\noindent The proof is sketched in Section~\ref{sec:asymptotics},
refer to Figure~\ref{fig:NumSim}, bottom right, for the graph
of~\eqref{eq:2}.

In particular, note that a stationary shock is present at $x=0$ in the
time asymptotic limit. Outside this shock, the solution is classical.

\section{Theorem \ref{th:example} -- Sketch of the Proof}
\label{sec:proof}

Note first that $H$, as defined in~\eqref{eq:hamiltonian},
satisfies~\eqref{eq:smoothness}, \eqref{eq:heterogeneity} with $X = 1$
and~\eqref{eq:convexity}. % With this flow, the conservation
% law~\eqref{eq:cl} reduces to the inviscid Burger equation with
% source term $-g'$. We attach to it the initial
% datum~\eqref{eq:initial} which evolves into a rarefaction in the
% homogeneous case. % \red{The
% proof of Theorem \ref{th:example} relies on the Cauchy problem for
% \eqref{eq:HS}, which reduces to, with $g$ as in
% \eqref{eq:hamiltonian}:}
% \begin{equation}
%   \label{eq:system}
%   \left\{
%     \begin{aligned}
%       \dot q & = p \\
%       \dot p & = -g'(q), \quad \mbox{with } g \mbox{ as in~\eqref{eq:hamiltonian}} \\
%       q(0) & = q_o \\
%       p(0) & = p_o.
%     \end{aligned}
%   \right.
% \end{equation}

% \subsection{Preliminary lemmas}
% \label{ssec:lemmas}

\begin{figure}[!h]
  \begin{center}
    \includegraphics[scale=0.66]{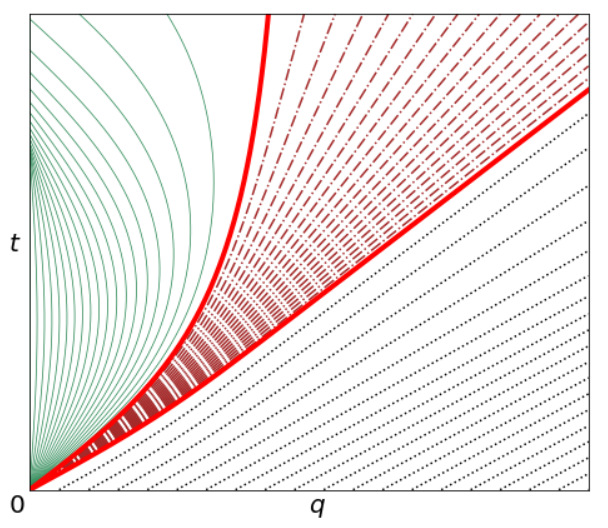}
  \end{center}
  \caption{On the horizontal axis, the $q$ component of solutions
    to~\eqref{eq:HS}--\eqref{eq:hamiltonian}, while time $t$ is on the
    vertical axis. Dotted curves are those considered in
    Lemma~\ref{lmm:ce1} and in Lemma~\ref{lmm:ce2}; dashed-dotted
    curves are those considered in \textit{(i)} of
    Lemma~\ref{lmm:ce3}; solid curves refer to \textit{(iii)} of
    Lemma~\ref{lmm:ce3}. The two thicker curves depict solutions
    corresponding to the initial data $(0, \sqrt{2})$ and $(0, 2)$.}
  \label{fig:curves}
\end{figure}

Here follows a sequence of lemmas describing the behaviors of the
solutions to~\eqref{eq:HS}--~\eqref{eq:hamiltonian}. The proofs,
relying on standard ODE arguments, can be found in~\cite{CPS2023}. We
refer to Figure~\ref{fig:curves} for an illustration of these
behaviors.

\begin{lemma}
  \label{lmm:ce1}
  {\rm{\cite[Lemma 5.11]{CPS2023}}} Let $H$ be as
  in~\eqref{eq:hamiltonian} and $u_o$ be as in~\eqref{eq:initial}.
  Fix $q_o \geq 0$. Denote by $(q, p)$ the solution
  to~\eqref{eq:HS}--\eqref{eq:hamiltonian} with initial datum
  $\left(q_o, u_o(q_o+)\right) = (q_o, 2)$. Then, $q$ is increasing
  on $[0, +\infty\mathclose[$ and $q(t) \limit{t}{+\infty} +\infty$.
\end{lemma}

\begin{lemma}
  \label{lmm:ce2}
  {\rm{\cite[Lemma 5.12]{CPS2023}}} Let $H$ be as
  in~\eqref{eq:hamiltonian} and $u_o$ be as in~\eqref{eq:initial}.
  Fix $0 \leq q_o < \widetilde{q_o}$ and denote by $(q, p)$,
  respectively $(\widetilde{q}, \widetilde{p})$, the global solution
  to~\eqref{eq:HS}--\eqref{eq:hamiltonian} with initial datum
  $\left(q_o, u_o(q_o+)\right) = (q_o, 2)$, respectively
  $\left(\widetilde{q_o}, u_o (\widetilde{q_o}+)\right) =
  (\widetilde{q_o}, 2)$. Then, $q(t) < \widetilde{q}(t)$, for all
  $t \geq 0$.
\end{lemma}

\begin{lemma}
  \label{lmm:ce3}
  {\rm{\cite[Lemmas 5.13~--~5.15]{CPS2023}}} Let $H$ be as
  in~\eqref{eq:hamiltonian} and $u_o$ be as
  in~\eqref{eq:initial}. Fix $p_o \in \mathopen]0,
  2\mathclose[$. Denote by $(q, p)$ the global solution
  to~\eqref{eq:HS}--\eqref{eq:hamiltonian} with initial datum $(0, p_o)$. \\[5pt]
  (i) If $p_o \in \mathopen]\sqrt{2}, 2\mathclose[$, then $q$ is
  increasing on
  $[0, +\infty\mathclose[$ and $q(t) \limit{t}{+\infty} +\infty$.\\[5pt]
  (ii) If $p_o = \sqrt{2}$, then $q$ is increasing on
  $[0, +\infty\mathclose[$ and
  $q(t) \limit{t}{+\infty} 1$. \\[5pt]
  % and $q$ is concave. \\[5pt]
  (iii) If $p_o \in \mathopen]0, \sqrt{2}\mathclose[$, then $q$ is
  periodic and the map
  \begin{equation}
    \label{eq:period}
    \fonction{\mathcal{T}}{\mathopen]0, \sqrt2\mathclose[}
    {\mathclose]0,+\infty\mathclose[}{p_o}{ \text{the smallest period of $q$}}
  \end{equation}
  is increasing, continuous and
  $\lim_{p_o \to \sqrt2} \mathcal{T} (p_o) = +\infty$.
  Moreover, $q (t)>0$ for
  $t \in \mathopen]0, \mathcal{T} (p_o)/2\mathclose[$ and $q (t)<0$
  for
  $t \in \mathopen]\mathcal{T} (p_o)/2, \mathcal{T} (p_o)\mathclose[$.
\end{lemma}

\begin{lemma}
  \label{lmm:ce4}
  {\rm{\cite[Lemma 5.16]{CPS2023}}} Let $H$ be as
  in~\eqref{eq:hamiltonian} and $u_o$ be as in~\eqref{eq:initial}.
  Fix $0 < p_o < \widetilde{p_o} < 2$ and denote by $(q, p)$,
  respectively $(\widetilde{q}, \widetilde{p})$, the global solution
  to~\eqref{eq:HS}--\eqref{eq:hamiltonian} with initial datum
  $(0, p_o)$, respectively $(0, \widetilde{p_o})$. Then,
  \[
    \begin{aligned}
      (i) \;\; & p_o \in ]0, \sqrt{2}[ \implies \forall t \in
      \mathopen]0, \mathcal{T} (p_o)/2], \;
      q(t) < \widetilde{q} (t) \\
      (ii) \;\; & p_o \in [\sqrt{2}, 2[ \implies \forall t \in
      ]0,+\infty[, \; q(t) < \widetilde{q} (t).
    \end{aligned}
  \]
\end{lemma}

% \subsection{Construction of the entropy solution by the method of
% characteristics}
% \label{ssec:EntropySol}

Recall that
\eqref{eq:smoothness}-\eqref{eq:heterogeneity}-\eqref{eq:convexity}
ensure that for all $(q_o, p_o) \in \R^2$, calling $(q,p)$ the
solution to~\eqref{eq:HS}--\eqref{eq:hamiltonian} with datum
$(q_o,p_o)$ at time $0$, the map
\begin{equation}
  \label{eq:flow}
  \fonction{\flow}{\R^3}{\R^2}{(t, q_o, p_o)}{\left(q(t), p(t)\right)}
\end{equation}
and its two projections $\flow_q, \flow_p$ are of class $\Ck{2}$, see
\cite[Lemma 5.2]{CPS2023}.

\begin{lemma}
  \label{lmm:ce5}
  {\rm{\cite[Lemma 5.17]{CPS2023}}} Let $H$ be as
  in~\eqref{eq:hamiltonian} and $u_o$ be as
  in~\eqref{eq:initial}. Then, there exists a unique map
  \begin{equation}
    \label{eq:mapDelta}
    \fonction{\Delta}{\mathopen]0, +\infty\mathclose[^2}{\left([0, +\infty\mathclose[ \times \{2\}\right) \cup
      \left(\{0\} \times \mathopen]0,2]\right)}{(t, x)}{(q_o, p_o)}
  \end{equation}
  such that
  \begin{equation}
    \label{eq:mapDelta2}
    \flow_q(t, q_o, p_o) = x
    \quad \text{and} \quad
    \forall \, s \in \mathopen]0, t\mathclose[, \; \flow_q(s, q_o, p_o) > 0 \,.
  \end{equation}
  Moreover,\\[5pt]
  (i) $\Delta$ is continuous.\\[5pt]
  (ii) $\Delta$ is monotone, in the sense that setting
  $\Delta (t_o,x_o) = (0,p_o)$ and $\Delta (t_o, x_o') = (0,p_o')$, if
  $0 < x_o < x_o'$, then $p_o < p_o'$.
  \\[5pt]
  (iii) For all $x \in \mathopen]0, +\infty\mathclose[$,
  $\lim_{t\to 0+} \Delta (t,x) = (x, 2)$.
\end{lemma}

\noindent The construction of the solution to~\eqref{eq:cl}--\eqref{eq:initial}
now follows from the next proposition.

\begin{proposition}
  \label{prop:ce6}
  {\rm{\cite[Proposition 5.18]{CPS2023}}} Let $H$ be as
  in~\eqref{eq:hamiltonian} and $u_o$ as
  in~\eqref{eq:initial}. Then,
  \begin{equation}
    \label{eq:entropysol}
    \fonction{u}{\mathopen]0, +\infty\mathclose[ \times (\R \setminus\{0\})}{\R}{(t,x)}{\left\{
        \begin{array}{ccl}
          \flow_p\left(t, \Delta(t, x)\right) & \text{if} & x > 0 \\[5pt]
          -\flow_p\left(t, \Delta(t, -x)\right) & \text{if} & x < 0.
        \end{array}
      \right.}
  \end{equation}
  is the solution to \eqref{eq:cl} in the sense of Definition
  \ref{def:EntropySolution} and it is a classical solution
  outside $x=0$.
\end{proposition}

We refer to~\cite[\S~5.4]{CPS2023} for the concluding remarks proving
the properties of $u$.

\section{Proof of Theorem~\ref{th:asymptotics}}
\label{sec:asymptotics}

Since $u$ is odd in the space variable, we consider only the case
$x>0$.  Using the notation~\eqref{eq:flow}, introduce the maps
\begin{equation}
  \label{eq:1}
  q^\flat (t) = \flow_q (t,0,\sqrt2)
  \quad\mbox{ and } \quad
  q^\sharp (t) = \flow_q (t,0,2)
\end{equation}
whose graphs are thicker in Figure~\ref{fig:curves}. By
Lemma~\ref{lmm:ce1} and Lemma~\ref{lmm:ce3},
$\lim_{t\to+\infty} q^\flat (t) = 1$ and
$\lim_{t\to+\infty} q^\sharp (t) = +\infty$.

Assume first $x \in \mathopen]0,1\mathclose[$. Then, if
$t > (q^\flat)^{-1} (x)$, by~(ii) in Lemma~\ref{lmm:ce3} and~(ii) in
Lemma~\ref{lmm:ce5}, $\Delta (t,x) = (0,p_o (t,x))$ for a suitable
$p_o(t,x) \in \mathopen]0, \sqrt 2 \mathclose[$. Moreover,
by~\eqref{eq:mapDelta2} in Lemma~\ref{lmm:ce5} and~(iii) in
Lemma~\ref{lmm:ce3}, $t < \mathcal{T} (p_o (t,x))/2$. Hence, as
$t \to +\infty$, also $\mathcal{T} (p_o (t,x)) \to +\infty$. By~(iii)
in Lemma~\ref{lmm:ce3}, $\lim_{t\to+\infty} p_o (t,x) = \sqrt2$.

Assume now that $x \geq 1$. If $t > (q^\sharp)^{-1} (x)$, by
Lemma~\ref{lmm:ce1}, $\Delta (t,x) = (0,p_o (t,x))$ with
$p_o (t,x) \in [\sqrt2,2\mathclose[$. By~(ii) in
Lemma~\ref{lmm:ce4} and~(i) in Lemma~\ref{lmm:ce3}, we have
$\lim_{t\to+\infty} p_o (t,x) = \sqrt2$.

Recall that along solutions to~\eqref{eq:HS}--\eqref{eq:hamiltonian},
$H$ is conserved, that is
\begin{equation}
  \label{eq:3}
  \frac{u(t,x)^2}{2} + g (x) = \frac{p_o (t,x)^2}{2} + g (0) \,,
\end{equation}
so that
\begin{equation}
  \label{eq:5}
\lim_{t\to+\infty} u (t,x)^2 = 2 (1-g (x)) \geq 0 \,,
\end{equation}
proving~\eqref{eq:2} for $x \geq 1$.

For $x \in [0,1\mathclose[$, we proved above that, for $t$ large,
$p_ o(t,x) \in \mathopen]0, \sqrt2\mathclose[$. This, together
with~\eqref{eq:3} ensures that
\begin{equation}
  \label{eq:4}
  u (t,x) < \sqrt{2 (1-g(x))} \,.
\end{equation}
By~\eqref{eq:cl}, understood in its classical sense thanks
to Proposition~\ref{prop:ce6}, and by~\eqref{eq:3}
\begin{displaymath}
  \partial_t u (t,x)
  =
  - \partial_x \left(H \left(x,u(t,x)\right)\right)
  =
  -\partial_x \left(\frac{p_o (t,x)^2}2\right)
  =
  -p_o (t,x) \;\partial_x p_o (t,x) \,.
\end{displaymath}
Indeed, note first that by~\eqref{eq:3}, $p_o (t,x) > 0$ and the
Implicit Function Theorem, $p_o$ is sufficiently regular. Item~(ii) in
Lemma~\ref{lmm:ce5} implies that $\partial_x p_o (t,x) \geq 0$, so
that $\partial_t u (t,x) \leq 0$. Hence, \eqref{eq:5} and~\eqref{eq:4}
allow to complete the proof of~\eqref{eq:2}.

\bigskip

As a final remark, we note that as $t\to+\infty$, $u (t)$ converges
uniformly on compact subsets of $\R \setminus \{0\}$,
thanks to pointwise convergence and monotonicity, by Dini's
Theorem~\cite[Theorem~7.13]{zbMATH03539473}. A further technical argument,
based on~(ii) of Theorem~\ref{th:example}, allows to extend the local uniform
convergence on all of $\R$.

\begin{acknowledgement}
  The first author was partly supported by the GNAMPA~2022 project
  \emph{Evolution Equations, Well Posedness, Control and
    Applications}.  This research was funded, in whole or in part, by
  l'Agence Nationale de la Recherche (ANR), project
  ANR-22-CE40-0010. For the purpose of open access, the second author
  has applied a CC-BY public copyright license to any Author Accepted
  Manuscript (AAM) version arising from this submission.
\end{acknowledgement}

% \newpage

\bibliographystyle{abbrv}

\bibliography{references_revised}

\begin{thebibliography}{10}

\bibitem{AdimurthiGhoshal2020}
{Adimurthi} and S.~S. Ghoshal.
\newblock Exact and optimal controllability for scalar conservation laws with
  discontinuous flux, 2020.

\bibitem{MR4188826}
F.~Ancona and M.~T. Chiri.
\newblock Attainable profiles for conservation laws with flux function
  spatially discontinuous at a single point.
\newblock {\em ESAIM Control Optim. Calc. Var.}, 26:Paper No. 124, 33, 2020.

\bibitem{MR1616586}
F.~Ancona and A.~Marson.
\newblock On the attainable set for scalar nonlinear conservation laws with
  boundary control.
\newblock {\em SIAM J. Control Optim.}, 36(1):290--312, 1998.

\bibitem{MR1722801}
E.~N. Barron, P.~Cannarsa, R.~Jensen, and C.~Sinestrari.
\newblock Regularity of {H}amilton-{J}acobi equations when forward is backward.
\newblock {\em Indiana Univ. Math. J.}, 48(2):385--409, 1999.

\bibitem{CP2020}
R.~M. Colombo and V.~Perrollaz.
\newblock Initial data identification in conservation laws and
  {H}amilton--{J}acobi equations.
\newblock {\em J. Math. Pures et Appl.}, \textbf{138}:{1--27}, 2020.

\bibitem{CPS2022}
R.~M. Colombo, V.~Perrollaz, and A.~Sylla.
\newblock Conservation laws and {H}amilton--{J}acobi equations with space
  inhomogeneity.
\newblock {\em Preprint}, 2022.
\newblock https://hal.archives-ouvertes.fr/hal-03873174.

\bibitem{CPS2023}
R.~M. Colombo, V.~Perrollaz, and A.~Sylla.
\newblock {I}nitial {D}ata {I}dentification in {S}pace {D}ependent
  {C}onservation {L}aws and {H}amilton-{J}acobi {E}quations.
\newblock {\em Preprint}, 2023.
\newblock https://hal.science/hal-04062783.

\bibitem{Dafermos1977}
C.~M. Dafermos.
\newblock Generalized characteristics and the structure of solutions of
  hyperbolic conservation laws.
\newblock {\em Indiana University Mathematics Journal},
  \textbf{26}(6):1097--1119, 1977.

\bibitem{MR4503821}
C.~Esteve-Yag\"{u}e and E.~Zuazua.
\newblock Reachable set for {H}amilton-{J}acobi equations with non-smooth
  {H}amiltonian and scalar conservation laws.
\newblock {\em Nonlinear Anal.}, 227:Paper No. 113167, 2023.

\bibitem{MR3643881}
L.~Gosse and E.~Zuazua.
\newblock Filtered gradient algorithms for inverse design problems of
  one-dimensional {B}urgers equation.
\newblock In {\em Innovative algorithms and analysis}, volume~16 of {\em
  Springer INdAM Ser.}, pages 197--227. Springer, Cham, 2017.

\bibitem{Kruzhkov1970}
S.~N. Kruzhkov.
\newblock First order quasilinear equations with several independent variables.
\newblock {\em Mathematics of the USSR-Sbornik}, \textbf{81}(123):228--255,
  1970.

\bibitem{LZ2021}
T.~Liard and E.~Zuazua.
\newblock Initial data identification for the one-dimensional {B}urgers
  equation.
\newblock {\em IEEE Transactions on Automatic Control}, 2021.

\bibitem{zbMATH03539473}
W.~Rudin.
\newblock Principles of mathematical analysis. 3rd ed.
\newblock International {Series} in {Pure} and {Applied} {Mathematics}.
  {D{\"u}sseldorf} etc.: {McGraw}-{Hill} {Book} {Company}. {X}, p.342, 1976.

\end{thebibliography}

\end{document}